\documentclass[a4paper,11pt]{amsart}
\usepackage[latin1]{inputenc}
\usepackage{amsmath,amsfonts,amssymb}
\def\R{\mathbb{R}}
\def\C{\mathbb{C}}
\def\N{\mathbb{N}}

\def\B{\mathbb{B}}

\def\rrr{\mathbb{R}}

\def\bbb{\mathbb{B}}

\def\sp{\vspace*{0.2cm}\\}

\def\dsp{\displaystyle}

\def\Im{\mathrm{Im}\,}

\def\Im{\mathrm{Im}\,}

\def\tsp{{\,}^t\!}
\def\vp{\varphi}
\def\ve{\varepsilon}
\def\D{\Delta}

\def\dist{\mathrm{dist}\,}

\theoremstyle{plain}
\newtheorem{theo}{Theorem}
\newtheorem{cor}{Corollary}
\newtheorem{prop}{Proposition}
\newtheorem{lm}{Lemma}
\newtheorem{defi}{Definition}
\newenvironment{dem}{\noindent{\bf{Proof}} \\}{\hfill $\square$ \newline}

\theoremstyle{remark}
\newtheorem{rem}{Remark}

\title[Regularity and estimates near the boundary for $J$-holomorphic discs]{Regularity and estimates for $J$-holomorphic discs attached to a maximal totally real submanifold}
\author{Léa Blanc-Centi}
\address{Léa Blanc-Centi, Universités de Marseille, Université de Provence, L.A.T.P., 39 rue Joliot-Curie, 13453 Marseille Cedex 13, FRANCE}
\email{lea@cmi.univ-mrs.fr}

\begin{document}
\begin{abstract}
We prove that pseudo-holomorphic discs attached to a maximal totally real submanifold inherit their regularity to the boundary from the regularity of the submanifold and of the almost complex structure. The proof is based on the computation of an explicit lower bound for the Kobayashi metric in almost complex manifolds, which also yields explicit estimates of Hölderian norms of such discs.
\end{abstract}

\maketitle
\section{Introduction}
Recent progress in symplectic geometry, and in particular the essential paper of M. Gromov \cite{Gr}, have strengthened the interest for almost complex manifolds and pseudo-holomorphic curves.
Pseudo-holomorphic discs form a natural invariant family for manifolds with boundary under the action of biholomorphisms, and more generally of CR maps. 
A {\em pseudo-holomorphic disc} in an almost complex manifold $M$ is a continuous map $h$ from the closed unit disc $\overline{\Delta}\subset\C$ to $M$, 
pseudo-holomorphic in $\D$. The first question is the existence of pseudo-holomorphic discs, and was solved by A. Nijenhuis et W. Woolf \cite{NW}. Indeed, they considered pseudo-holomorphic discs as the solutions of non-linear elliptic operators, and proved that at any point, in any direction, there is a small pseudo-holomorphic disc. This argument also shows that the regularity of the disc in $\D$ depends on the regularity of the almost complex structure, and it yields \textsl{a priori} estimates \cite{S}.\\ 

The aim of this paper is to study boundary properties of pseudo-holomorphic discs. We will consider discs $h$ {\em attached} to a submanifold $E$, that is, such that $h(\partial\D)\subset E$. 
The boundary properties of such discs are strongly related to the geometry of the manifold to which they are attached.
 
In the complex case, discs attached to a maximal totally real submanifold have very interesting regularity properties \cite{Chirka,CCS,IV}, which come from a reflexion principle.
In the almost complex situation, we know by \cite{CGS} that the discs are smooth up to the boundary if the submanifold $E$ and the almost complex structure $J$ are smooth. We prove a quantitative version of this result.
In fact, when $E$ and $J$ are only supposed to have Hölderian regularity, we show that the discs inherit from the minimum regularity between the one of $E$ and the one of $J$. We also give explicit estimates of Hölderian norms of the discs. More precisely:

\begin{theo}\label{p4.7-CGS}
Let $(M,J)$ be a $\mathcal{C}^r$-almost complex manifold (where $r\ge 1$ is not an integer) and $E\subset M$ be a maximal totally real submanifold.\\ 
Then, every continuous map $h$ from $\Delta^+\cup ]-1;1[$ to $M$, $J$-holomorphic on $\Delta^+$ and attached by the diameter to $E$, is locally of class $\mathcal{C}^{r}$ in  $\Delta^+\cup
  ]-1;1[$.
Moreover, for any compact subset $K$ in $\Delta^+\cup
  ]-1;1[$, 
$$||h||_{\mathcal{C}^{r}(K)}\le c(r,K)||h||_\infty\left(1+\frac{c(K)}{\sqrt{\lambda^J_E}}\right).$$
\end{theo}

\noindent Here, $\lambda_E^J$ denotes the smallest eigenvalue of the Levi form of a function of the type ``squared distance to $E$''. In view of its geometric interpretation, we call $\lambda_E^J$ the \textit{minimal $J$-curvature of $E$}. Note that the estimate given by L. Lempert in \cite{L1} for stationnary discs in strongly convex domains of $\C^n$ also used the curvature of the domain.

The proof is based on the computation of an explicit lower bound for the Kobayashi-Royden infinitesimal pseudometric in almost complex manifolds (Proposition \ref{prop:minorationK'}).
 This lower bound is obtained by constructing $J$-plurisubharmonic maps, and yields an estimate of the size of pseudo-holomorphic discs (Proposition \ref{prop:tailledisques}).\\

Finally, we apply Theorem \ref{p4.7-CGS} to the study of a pseudo-holomorphic map along the edge of a wedge. The method, originally introduced by E. Bishop \cite{Bishop}, consists in relating the behavior in the interior and the behavior at the boundary by means of analytic discs. The idea is to fill the interior with analytic discs ``glued'' to the edge. We get:

\begin{cor}\label{cor:wedge}
Let $(M,J)$ and $(M',J')$ be respectively a $\mathcal{C}^{r}$ and a $\mathcal{C}^{r'}$-almost complex manifold (where $r,r'\ge 1$ are not integers), and $\Omega\subset M$ be a domain. We suppose that $N\subset \Omega$ is a maximal totally real submanifold of $(M,J)$ of class $\mathcal{C}^r$, and that $N'$ is a maximal totally real submanifold of $(M',J')$ of class $\mathcal{C}^{r'}$. Keep  
$s=\mathrm{min}\,(r-1,r')$.\\
Then every pseudo-holomorphic map $F:W(\Omega,N)\to (M',J')$, continuous on $W(\Omega,N)\cup N$, and such that $F(N)\subset
N'$, is locally of class $\mathcal{C}^s$ and verifies for any compact set $K$ in $W(\Omega,N)\cup N$: 
$$||F||_{\mathcal{C}^{s}(K)}\le c(s,K)||F||_\infty\left(1+\frac{c(K)}{\sqrt{\lambda_{N'}^{J'}}}\right).$$
\end{cor}

This paper is organized as follows. In the second section, we recall some basic notions of almost complex geometry. The third section is devoted to the obtaining of an explicit lower bound of the Kobayashi metric; we also give an explicit estimate of the size of pseudo-holomorphic discs. In Section 4, we establish Theorem \ref{p4.7-CGS} and Corollary \ref{cor:wedge}.



\section{Preliminaries}
\subsection{Almost complex manifolds and pseudo-holomorphy}
Let us recall some definitions.

\begin{defi}
An {\em almost complex structure} on a smooth (real) manifold $M^{2n}$ is a $(1,1)$-tensor $J$, that is, a section from $M$ to $End(TM)$, such that $J^2=-Id$. If $J$ is of class $\mathcal{C}^r$, we say that $(M,J)$ is a {\em $\mathcal{C}^r$-almost complex manifold}.
\end{defi}

The first example of an almost complex manifold is the space $\R^{2n}$ equipped with the {\em standard complex structure} $J_{st}$ given at any point, in canonical coordinates, by the matrix $J_{st}=\left(\begin{array}{cc}0& -I_n\\ I_n&0\end{array}\right)$
where $I_n$ is the $n$-sized identity matrix. 

We will also consider the following situation:

\begin{defi}
A submanifold $N$ in an almost complex manifold $(M^{2n},J)$ is {\em maximal totally real} if $N$ is of (real) dimension $n$ and if $TN\cap J\,TN=\{0\}$.
\end{defi}

\begin{defi}
A $\mathcal{C}^1$ regular map $F:(M,J)\to(M',J')$ between almost complex manifolds is said to be {\em $(J,J')$-holomorphic} if $J'\circ dF=dF\circ J$.
\end{defi}

If $(M,J)$ is the unit disc of $\C$ (that is, $\D\subset\R^2$ equipped with the standard complex structure), we say that $F$ is a {\em $J'$-holomorphic}, or {\em pseudo-holomorphic, disc}.
Generically, given $(M,J)$, there does not always exist pseudo-holomorphic maps from $(M,J)$ to $(M',J')$, even if $(M',J')=(\Delta,J_{st})$. Nevertheless, when $(M,J)=(\Delta,J_{st})$, the condition of pseudo-holomorphy can be translated into a quasi-linear elliptic system of first order that admits non-trivial solutions \cite{NW}. Hence there exist infinitely many pseudo-holomorphic discs $h$ valued in a fixed almost complex manifold. Moreover, one can prescribe $h(0)$ and $dh_0(\partial/\partial\, x)$.

As in the complex situation, maps verifying the equation of pseudo-holomorphy inherit their smoothness from the smoothness of the almost complex structures: if $J$ and $J'$ are of class $\mathcal{C}^r$, then every $(J,J')$-holomorphic map  is of class $\mathcal{C}^{r+1}$.

\subsection{The Levi form}
Let $(M,J)$ be an almost complex manifold. The {\em Levi form} $\mathcal{L}^Ju$ of any $\mathcal{C}^2$-regular function $u:M\to\R$ is defined by
$$\mathcal{L}^Ju(X)=d(d^c_Ju)(X,JX),\quad\mathrm{where}\quad d_J^cu(X)=-du(JX)$$
for all $X\in\ TM$.
In local coordinates, $\dsp d_J^cu=-\sum_{i,j}\frac{\partial u}{\partial
  x_i}\,J_{i,j}dx_j$ and
$$d(d_J^cu)=-\sum_{i,j,k}\frac{\partial^2u}{\partial x_i\partial
  x_k}\,J_{i,j}dx_k\wedge dx_i-\sum_{i,j,k}\frac{\partial u}{\partial
  x_i}\,\frac{\partial J_{i,j}}{\partial x_k}\,dx_k\wedge dx_j\,.$$
Since $dx_k\wedge
dx_j(X,JX)=x_k(Jx)_j-(Jx)_kx_j=\sum_lx_kx_lJ_{j,l}-\sum_lx_jx_lJ_{k,l}$,
we easily get the following expression:

\begin{eqnarray}\label{expressionLevi}
\mathcal{L}^Ju(X)=\tsp XDX+\tsp(JX)D(JX)+\tsp X(A-\tsp A)JX,
\end{eqnarray}
where 
$$D=\left(\frac{\partial^2u}{\partial x_i\partial
    x_j}\right)_{1\le i,j\le 2n}\quad \mathrm{and}\quad
A=\left(\sum_i\frac{\partial u}{\partial x_i}\frac{\partial
    J_{i,j}}{\partial x_k}\right)_{1\le j,k\le 2n}.$$

If $J_0$ is any constant almost complex structure, this gives, by setting $J=J_0+H$:
\begin{equation}\label{etoile}
\mathcal{L}^Ju(X)=\mathcal{L}^{J_0}u(X)+2{\,}\tsp(HX)D(J_0X)+\tsp(HX)D(HX)+\tsp
X(A-\tsp A)(J_0+H)X.
\end{equation}

Note that if $\phi:x\mapsto P^{-1}x$ is a linear change of variables in $\R^{2n}$, (\ref{expressionLevi}) shows that for   
$\tilde{J}=\phi_*J=P^{-1}J$ and $\tilde{u}=u\circ\phi^{-1}$, 
$$\mathcal{L}^{\tilde{J}}\tilde{u}_x(X)=\tsp
X\tilde{D}(x)X+\tsp(\tilde{J}(x)X)\tilde{D}(x)(\tilde{J}(x)X)+\tsp X\tilde{A}(x)\tilde{J}(x)X-\tsp(\tilde{J}(x)X)\tilde{A}(x)X$$ 
where
$\displaystyle\tilde{D}(x)=\left(\frac{\partial^2\tilde{u}}{\partial x_i\partial x_j}(x)\right)_{i,j}=\tsp PD(Px)P$ and $\displaystyle\tilde{A}(x)=\left(\sum_i\frac{\partial\tilde{u}}{\partial x_i}\frac{\partial\tilde{J}_{i,j}}{\partial x_k}\right)_{j,k}=\tsp PA(Px)P$, and so
\begin{equation}\label{chgtcartelin}
\mathcal{L}^{\tilde{J}}_x\tilde{u}(X)=\mathcal{L}^J_{Px}u(PX).
\end{equation}

\subsection{Small perturbations of the standard structure}\label{s:SmallPerturbation}
Let $(M,J)$ be an almost complex manifold of dimension $2n$. We can assume locally that $M$ is an open subset in $\R^{2n}$, and that $J$ takes its values in $\mathcal{M}_{2n}(\R)$: we write $J=\left(\begin{array}{cc}A&B\\C&D\end{array}\right)$ where $A,B,C,D$ are $n$-sized blocks.
If $J$ is sufficiently near $J_{st}$, then 
\begin{eqnarray}\label{petiteperturbation}
J^2=-I_{2n} \Longleftrightarrow J=\left(\begin{array}{cc}A&-(I_n+A^2)C^{-1}\\C&-CAC^{-1}\end{array}\right)\,.
\end{eqnarray}

The two following lemmas allow us to choose local charts with special properties. We give the proofs in the Appendix.

\begin{lm}\label{lm:chgtcarte}
Let $(M^{2n},J)$ be a $\mathcal{C}^r$-almost complex manifold, $E$ a maximal totally real $\mathcal{C}^{r'}$-submanifold (where $r,r'\ge 1$). Fix $\ve>0$ and $p\in M$.
Then there exist local coordinates $z:U\to\B^{2n}$ in a neighborhood of $p$ such that $||z_*J-J_{st}||_{\mathcal{C}^1}\le\varepsilon$. If one also assumes $p\in E$, the map $z$ is of class $\mathcal{C}^{\mathrm{min}(r,r')}$ and verifies:
\begin{itemize}
\item $z(U\cap E)=\R^n\cap \B^{2n}$;
\item $||{(z_*J)}_{(x^*,y^*)}-J_{st}||\le c\ \varepsilon\,||y^*||$ for some constant $c$;
\item the equation of $z_*J$-holomorphy is $\displaystyle\frac{\partial}{\partial\bar{\zeta}}+Q\times\frac{\partial}{\partial\zeta}=0$, where $Q$ is a $\mathcal{C}^{\mathrm{min}(r,r')-1}$ regular map such that $||Q(x^*,y^*)||\le c\ \varepsilon\,||y^*||$ for some constant $c$.
\end{itemize} 
\end{lm}

\noindent The corresponding basis of $(1,0)$-forms can then be written 
$dz^*+O(||y^*||)dy^*$.

We will say that an atlas of $M$ in a neighborhood of $E$ is \textit{$(\ve,E)$-tamed} if it consists of charts verifying the properties of Lemma \ref{lm:chgtcarte}. Since a bounded domain $D\subset M$ admits a finit covering $(z_i,U_i)_{1\le i\le s}$ by such charts, we set for every $\mathcal{C}^{s}$-regular map that takes its values in $\overline{D}$:
$$||f||_\infty=\underset{i}{\mathrm{Max}}\,||z_i\circ f_{|f^{-1}(U_i)}||_\infty\quad\mathrm{and}\quad||f||_{s}=\underset{i}{\mathrm{Max}}\,||z_i\circ f_{|f^{-1}(U_i)}||_{s}.$$

It is possible to lift the almost complex structure to the tangent bundle. Moreover, the properties listed in Lemma \ref{lm:chgtcarte} are preserved:

\begin{lm}\label{lm:relevement}
Let $(M,J)$ be a $\mathcal{C}^r$-almost complex manifold (where $r\ge 2$ is not an integer), and $E$ a maximal totally real manifold and $\ve>0$.
Assume that $p\in E$ and that $(z,U)$ is a $(\ve,E)$-tamed chart in a neighborhood of $p$. 
Set $M^c=TM$ and $E^c=TE$.\\
Then there exists an almost complex structure $J^c$ on $M^c$, $\mathcal{C}^{r-1}$ regular, that induces $J$ on $M$ and such that:\\
\textbullet\ $E^c$ is a maximal totally real submanifold of $(M^c,J^c)$;\\
\textbullet\ $z^c=(z,dz)$ is $(\ve,E^c)$-tamed in a neighborhood of $(p,0)\in E^c$.
\end{lm}

\section{An explicit lower bound for the Kobayashi metric}
\subsection{Estimates of the Levi form for some classical functions}
In view of Lemma \ref{lm:chgtcarte}, we may mainly consider small perturbations of the standard structure. In order to measure the perturbation, we introduce the following notations:
$$||J_p||_0=\mathrm{Max}\{||J_pX||/\ X\in\mathcal{M}_{2n,1}(\R),\ ||X||=1\}$$
$$||J_p||_1=||J_p||_0+\left(\sum_{i=1}^{2n}\left|\left|\left(\frac{\partial
          J_{i,j}}{\partial x_k}\right)_{j,k}(p)\right|\right|_0^2\right)^{1/2}$$
where $||X||=(\sum_{i=1}^{2n}x_i^2)^{1/2}$ denotes the Euclidian norm. Moreover, for some relatively compact subset $D$ in $M$, we set
$$||J||_{\mathcal{C}^1(\bar{D})}=\mathrm{Max}\{||J_p||_1/\ p\in\bar{D}\}.$$
This defines a norm on $\mathcal{C}^r(M,\mathcal{M}_{2n}(\R))$, and provides the classical $\mathcal{C}^1$-topology.

Set $J=J_{st}+H$, and let $u:M\to\R$ be of class $\mathcal{C}^2$. Then, we get by (\ref{expressionLevi}) that for all $p\in M$ and $X\in T_pM$:
$$\mathcal{L}_p^Ju(X)\ge\mathcal{L}^{J_0}_pu(X)-2\rho(p)||H||_0||J_{st}||_0||X||^2+\mu(p)||HX||^2||A||_0||J||_0||X||^2.$$
Since $A=\sum_i\frac{\partial u}{\partial x_i}\left(\frac{\partial
    H_{i,j}}{\partial x_k}\right)_{j,k}$, we have
$||A||_0\le\sqrt{\sum_i\left|\frac{\partial u}{\partial
      x_i}\right|^2}\times\sqrt{\sum_i\left|\left|\left(\frac{\partial
          H_{i,j}}{\partial x_k}\right)_{j,k}\right|\right|_0^2}$, and thus $||A(p)||_0\le ||\nabla u_p||\times ||H(p)||_1$.
Hence we get:
\begin{eqnarray}\label{minorationLevi}
\lefteqn{\mathcal{L}^J_pu(X)\ge\mathcal{L}^{J_{st}}_pu(X)-2\rho(p)||H(p)||_0||X||^2+\mu(p)||H(p)X||^2}\hspace{4cm}\nonumber\\
 &-2||\nabla u_p||(1+||H(p)||_0)||H(p)||_1\times||X||^2
\end{eqnarray}
where $\rho(p)=\mathrm{Max}|\lambda|$ and $\mu(p)=\mathrm{min}\,\lambda$ for $\lambda$ describing the eigenvalues of the real Hessian matrix $\left(\frac{\partial^2u}{\partial x_i\partial x_j}(p)\right)_{i,j}$.

We first search an explicit lower bound for the Levi form of each function $||z||$,\ $||z||^2$,\ and $\ln||z||$. 
We thus need to determine the eigenvalues of the corresponding Hessian matrices. We will use that for every non-zero column vector $V$, the matrix $\lambda(I-V\tsp V)$ has exactly two eigenvalues: 

\textbullet\ $\lambda(1-||V||^2)$, whose eigenspace is $\mathrm{Vect}\,V$;

\textbullet\ $\lambda$, whose eigenspace is $V^\perp=\{X/\ \tsp VX=0\}$.\\
\noindent Indeed, the eigenvalues of $\lambda(I-\tsp VV)$ can be written under the form $\lambda(1-\mu)$, where $\mu$ is any eigenvalue of the matrix $V^tV$.
Suppose $V\tsp VX=\mu X$: since $\tsp VX$ is scalar, then either $\mu=0$ (and $\tsp VX=0$), or $X$ is proportional to $V$. Conversely, if $X=V$, then $V\tsp VX=V\times(\tsp VV)=||V||^2X$.\\

\noindent We obtain by (\ref{minorationLevi}) the following lower estimates:
\begin{eqnarray*}
\ \mathcal{L}^J_p(||z||)(X)\!\!&\!\!\ge\!\!&\!\!||X||^2\left(\!\frac{1}{||p||}-\frac{2}{||p||}||H(p)||_0-2(1+||H(p)||_0)\times||H(p)||_1\!\right)\\
\ \mathcal{L}^J_p(\ln||z||)(X)\!\!&\!\!\ge\!\!&\!\!||X||^2\left(\!\!-\frac{2}{||p||^2}||H(p)||_0-\frac{1}{||p||^2}||H(p)||_0^2-\frac{2}{||p||}(1+||H(p)||_0)\!\times\!||H(p)||_1\!\!\right)\\
\ \mathcal{L}^J_p(||z||^2)(X)\!\!&\!\!\ge\!\!&\!\!||X||^2\left(4-4||H(p)||_0-4||p||(1+||H(p)||_0)\times||H(p)||_1\right)\,.
\end{eqnarray*}

\begin{rem}
Suppose that $p$ is in some domain $D\ni 0$ of the unit ball, and that $J(0)$ is the standard strcture. Then
$||J_p-J_{st}||\le||p||\times\underset{\bar{D}}{\mathrm{Max}}|||dJ_p|||$.
Let us define $||J||_{\bar{D}}=\underset{\bar{D}}{\mathrm{Max}}|||dJ_p|||+||J||_{\mathcal{C}^1(\bar{D})}$:
\begin{eqnarray*}
\mathcal{L}^J_p(||z||)(X)&\ge&||X||^2\left(\frac{1}{||p||}-2||H||_{\bar{D}}-2(1+||H||_{\bar{D}}\times||p||)\times||H||_{\bar{D}}\right)\\
\mathcal{L}^J_p(\ln||z||)(X)&\ge&-||X||^2\left(\frac{2}{||p||}||H||_{\bar{D}}+||H||_{\bar{D}}^2+\frac{2}{||p||}(1+||H||_{\bar{D}}||p||)\times||H||_{\bar{D}}\right)\\
\mathcal{L}^J_p(||z||^2)(X)&\ge&||X||^2\left(4-4||H||_{\bar{D}}\times||p||-4||p||(1+||H||_{\bar{D}}\times||p||)\times||H||_{\bar{D}}\right).
\end{eqnarray*}
\end{rem}

\subsection{Construction of strictly $J$-plurisubharmonic functions}
For $D$ being a bounded domain in $M$ and $u:\bar{D}\to\R$ a function of class $\mathcal{C}^2$, the continuous function 
$$(x,X)\mapsto\frac{\mathcal{L}^J_xu(X)}{||X||^2}$$
admits a minimum
on the compact set $\bar{D}\times S^{2n-1}$, which is the lower eigenvalue of the quadratic form $\mathcal{L}^Ju$ on
$\bar{D}$. We will denote it by $\lambda_0(D,J,u)$.

\begin{defi}
The function $u$ is said to be {\em $J$-plurisubharmonic} on $D$ when $\lambda_0(D,J,u)\ge 0$, that is, its Levi form is positive semi-definite. It is {\em strictly $J$-plurisubharmonic} on $D$ when $\lambda_0(D,J,u)> 0$, that is, its Levi form is positive definite.
\end{defi}
 
In fact, $u$ is strictly $J$-plurisubharmonic if and only if for all $J$-holomorphic disc $h$ in $M$ centered at $p$ and such that $(\partial h/\partial x)(0)=v\not=0$:
$$\Delta(u\circ h)_0=\mathcal{L}^J_pu(v)>0.$$

The method in order to construct strictly $J$-plurisubharmonic functions consists in perturbing some strictly $J_{st}$-plurisubharmonic functions. Here we work in local coordinates, and we denote by $J_{st}=\left(\begin{array}{cc}0&-I_n\\ I_n&0\end{array}\right)$ the standard complex structure. 

\begin{lm}\label{lm:epsilonm}
Let $D$ be a bounded domain in $\R^{2n}$ with an almost complex structure
$J$. Suppose that $\theta$ is a non-decreasing smooth function on $\R^+$, such that 
$\theta(x)=x$ for $x\le 1/3$ and $\theta(x)=1$ for $x\ge 2/3$.
If $D$ is bounded by $m$, we set $\varepsilon_m=\mathrm{min}\left(\frac{1}{32(1+m)},\frac{1}{32m(1+m)}\right)$. Then:
\begin{itemize}
\item For all $x\in \bar{D}$, $p\in \bar{D}$ and $X\in\R^{2n}$,
$$||J-J_{st}||_{\mathcal{C}^1(\bar{D})}\le\varepsilon_m\Longrightarrow\
    \frac{7}{2}||X||^2\le\mathcal{L}^J_x(||\cdot-p||^2)(X)\le \frac{9}{2}||X||^2.$$
\item There exists some constant $k$ depending only on $\theta$, such that for all $r>0$, $A>1$, $B\ge k$, $p\in D$, and for every almost complex structure $J$ verifying
$J(p)=J_{st}$ and $||J-J_{st}||_{\mathcal{C}^1(\bar{D})}\le\varepsilon_m$, the function defined by
$$x\overset{u}{\mapsto}\ln\left(\theta\left(\frac{1}{r^2}||x-p||^2\right)\right)+A||x-p||+B\frac{1}{r^2}||x-p||^2$$
is strictly $J$-plurisubharmonic on $D$.\sp
\item Assume that $||J-J_{st}||_{\mathcal{C}^1(\bar{D})}\le\varepsilon_m$ and let $w:\bar{D}\to\R$ be a $J$-plurisubharmonic function of class $\mathcal{C}^2$. Then, for all $\delta\le\frac{2}{9}\lambda_0(D,J,w)$ and $p\in D$, the map
$$x\mapsto w(x)-\delta||x-p||^2$$
is strictly $J$-plurisubharmonic on $D$.
\end{itemize}
\end{lm}

\begin{dem}
\noindent\textbullet\ Set $H=J-J_{st}$: then $||H||_{\mathcal{C}^1(\bar{D})}\le\varepsilon_m$. We obtain:
$$\mathcal{L}^J_x(||\cdot-p||^2)\ge||X||^2(4-4||H||_{\mathcal{C}^1(\bar{D})}-4||x-p||(1+||H||_{\mathcal{C}^1(\bar{D})})||H||_{\mathcal{C}^1(\bar{D})}),$$
and, in the same way,
$$\mathcal{L}^J_x(||\cdot-p||^2)\le||X||^2(4+4||H||_{\mathcal{C}^1(\bar{D})}+2||H||_{\mathcal{C}^1(\bar{D})}^2+(1+||H||_{\mathcal{C}^1(\bar{D})}))\times 4||x-p||\,||H||_{\mathcal{C}^1(\bar{D})}.$$
Since we can assume $\ve_m\le 1$, we get:
$$4-(16+16m)\varepsilon_m\le \frac{\mathcal{L}^J_x(||\cdot-p||^2)}{||X||^2}\le 4+(16+16m)\varepsilon_m.$$
\textbullet\ Pick $p\in D$ and set
$v:x\mapsto\ln\left(\theta\left(\left|\left|\frac{x-p}{r}\right|\right|\right)^2\right)$.
We have
$$||\nabla
v(x)||=\frac{\theta'}{\theta}\left(\frac{||x-p||^2}{r^2}\right)\times
2\frac{||x-p||}{r^2}$$
and, by setting $V=\tsp (x_1-x_1(p),\hdots,x_{2n}-x_{2n}(p))$, the matrix $\left(\frac{\partial^2v}{\partial x_i\partial
  x_j}(x)\right)_{i,j}$ is equal to $\alpha I_{2n}+\beta V^tV$ where
$$\alpha=\frac{\theta'}{\theta}\left(\frac{||x-p||^2}{r^2}\right)\times\frac{2}{r^2}\quad\mathrm{and}\quad\beta=\left(\frac{\theta''\theta-\theta'^2}{\theta^2}\right)\left(\frac{||x-p||^2}{r^2}\right)\times\frac{4}{r^4}.$$

\noindent This gives the two eigenvalues of the real Hessian of $v$, and hence a lower bound for its Levi form.\\
Since $v(x)=2(\ln||x-p||-\ln r)$ on $D\cap \B(p,r/\sqrt{3})$,
$$\mathcal{L}^J_xv(X)\ge -||X||^2\times 8\left(\frac{1}{||x-p||}+1\right)\,\varepsilon_m.$$
The map vanishes on $D\setminus \B(p,r)$, thus $\mathcal{L}^J_xv(X)=0$.\\
\ \\
On $D\cap \B(p,r)\setminus \B(p,r/\sqrt{3})$: 
$$\mathcal{L}^J_xv(X)\ge -||X||^2(\rho(x)(2+\varepsilon_m)+2\varepsilon_m(1+\varepsilon_m)||\nabla v||).$$
Let us set
\begin{eqnarray}\label{cstek}
k=4\times\mathrm{Max}\left(\underset{\frac{1}{3}\le
    x\le 1}{\mathrm{Sup}}\left|\frac{\theta'}{\theta}\right|,\underset{\frac{1}{3}\le
    x\le
    1}{\mathrm{Sup}}\left|\frac{\theta''\theta-\theta'^2}{\theta^2}\right|\right).
\end{eqnarray}
 Then, in $D\cap \B(p,r)\setminus \B(p,r/\sqrt{3})$, we have
$$\rho(x)\le
\frac{k}{4}\times\left(\frac{2}{r^2}+\frac{4||x-p||^2}{r^4}\right)\ \
\mathrm{and}\ \ ||\nabla v(x)||\le \frac{k}{4}\times\frac{2||x-p||}{r^2},$$
and hence
\begin{equation}
\mathcal{L}^J_xv(X)\ge-||X||^2\times\frac{k}{2r^2}\left((1+\frac{2}{r^2}||x-p||^2)(2+\varepsilon_m)+2\varepsilon_m(1+\varepsilon_m)||x-p||\right).
\end{equation}
Now we study the Levi form of $u:x\mapsto u(x)=v(x)+A||x-p||+B\left|\left|\frac{x-p}{r}\right|\right|^2$.
The choice of $\varepsilon_m$ implies that on $D\cap \B(p,r/\sqrt{3})$, 
$$\frac{\mathcal{L}^J_xu(X)}{||X||^2}\ge\frac{1}{2||x-p||}(A-1)$$
since $\varepsilon_m=\mathrm{min}\left(\frac{1}{32(1+m)},\frac{1}{32m(1+m)}\right)$.\vspace*{0.1cm}\\
One gets 
$$\frac{\mathcal{L}^J_xu(X)}{||X||^2}\ge\frac{1}{2||x-p||}A$$
on $D\setminus \B(p,r)$ and 
$$\frac{\mathcal{L}^J_xu(X)}{||X||^2}\ge\frac{1}{2||x-p||}A+\frac{7}{2r^2}(B-k)$$
on $D\cap \B(p,r)\setminus \B(p,r/\sqrt{3})$.

Consequently, for $A>1$ and $B\ge k$, the Levi form $\mathcal{L}_xu$ is positive definite at every $x\in D$.\\
\ \\
\textbullet\ Let $p\in D$ and $\tilde{w}:x\mapsto w(x)-\delta||x-p||^2$:
$$\forall x\in D,\
\mathcal{L}^J_x\tilde{w}(X)\ge\mathcal{L}^J_xw(X)-\frac{9}{2}\delta||X||^2
 \ge\left(\lambda_0(D,J,w)-\frac{9}{2}\delta\right)||X||^2.$$
\end{dem}

\subsection{Estimate of the Kobayashi metric}
Let $(M,J)$ be an almost complex manifold. 
For every $p\in D$ and $v$ a tangent vector at point $p$, we set
$$K_{(D,J)}(p,v)=\text{inf}\{\alpha>0/\ \exists h\in\mathcal{O}^J(\Delta,D)\ \text{with}\ h(0)=p\ \text{and}\ (\partial h/\partial x)(0)=v/\alpha\},$$
which is well-defined according to \cite{NW}. 
Most of the basic properties of the Kobayashi-Royden pseudo-metric in the complex situation are always true in the almost complex case, as the decreasing property under the action of pseudo-holomorphic maps. 
\\
\ \\
Let $D$ be a bounded domain in $\R^{2n}$ equipped with an almost complex structure $J$, and $u$ be any negative $\mathcal{C}^2$-regular function on $\bar{D}$,
strictly $J$-plurisubharmonic on $D$. We assume that $J(p)=J_{st}$ and
$||J-J_{st}||_{\mathcal{C}^1(\bar{D})}\le\varepsilon_m$ with $m$ a bound of $D$ and $\varepsilon_m$ given by Lemma \ref{lm:epsilonm}. By means of Lemma \ref{lm:epsilonm}, the proof of Proposition 4.4 in \cite{CGS} shows that there exists some constant $c_m=\sqrt{\frac{2}{9ke^{2m}}}$ depending only on $m$ such that 
$\displaystyle K_{(D,J)}(p,v)\ge
\sqrt{\frac{\delta}{eBe^{2Am}}}\,\frac{||v||}{\sqrt{|u(p)|}}$ for every
 $A>1$, $B\ge k$, $\delta<\frac{2}{9}\lambda_0(D,J,u)$. Hence
\begin{equation}\label{minorationK}
K_{(D,J)}(p,v)\ge c_m\sqrt{\lambda_0(D,J,u)}\times\frac{||v||}{\sqrt{|u(p)|}}
\end{equation}
where $k$ is defined by (\ref{cstek}). 
The next step is to remove the hypothesis $J(p)=J_{st}$, in order to get a uniform estimate:

\begin{prop}\label{prop:minorationK'}
Let $D$ be a domain in the unit ball $\B^{2n}$ equipped with an almost complex structure $J$. Assume that $u$ is any negative $\mathcal{C}^2$ regular function on $\bar{D}$, 
strictly $J$-plurisubharmonic on $D$.\\ 
There exist some universal constants $c'$ and $\varepsilon'$ such that if $||J-J_{st}||_{\mathcal{C}^1(\bar{D})}\le\varepsilon'$, then
$$\forall p\in D,\ \forall v\in\R^{2n},\ 
K_{(D,J)}(p,v)\ge
c'\,e^{-2t}\,\sqrt{\lambda_0(D,J,u)}\,\frac{||v||}{\sqrt{|u(p)|}}$$
as soon as $D$ is bounded by $t$.
\end{prop}

\begin{dem}
Pick $p\in D$. The matrix $P_p$ consisting of the vectors
$(e_1,\hdots,e_n,J_pe_{n+1},\hdots,J_pe_{2n})$, where $(e_1,\hdots,e_n,e_{n+1},\hdots,e_{2n})$ is the canonical basis of 
$\R^{2n}$, depends continuously on $J$ and $p$. 
So one can fix $\varepsilon'$ such that $||J-J_{st}||_{\mathcal{C}^1(\bar{D})}\le\varepsilon'$ implies that for all $p\in\B^{2n}$, the matrix $P_p$ is invertible, $||P_p^{-1}JP_p-J_{st}||_{\mathcal{C}^1(\bar{D})}\le\varepsilon_{m=2}$ and $||P_p^{-1}||_0,||P_p||_0\le 2$. 

If $D$ is bounded by $t\in]0;1]$, then $||P_p^{-1}JP_p-J_{st}||_{\mathcal{C}^1(\bar{D})}\le\varepsilon_{m=2t}$ since $\ve_m$ is a decreasing function of $m$.
For $\phi:x\mapsto P_p^{-1}x$, we get $\phi_*J(\phi(p))=J_{st}$ and $\phi(D)\subset 2t \B^{2n}$. Thus (\ref{minorationK}) applies to $\tilde{D}=\phi(D),\ \tilde{J}=\phi_*J,\ \tilde{u}=u\circ\phi^{-1}$ and gives 
$$K_{(D,J)}(p,v)\ge c_{m=2t}\sqrt{\lambda_0(\tilde{D},\tilde{J},\tilde{u})}\times\,\frac{||d(\phi_p)_p(v)||}{\sqrt{|u(p)|}}.$$
Since $||d(\phi_p)_p(v)||=||P_p^{-1}v||\ge||v||/||P_p||_0\ge||v||/2$ in view of the choice of $\varepsilon'$, we obtain
$$K_{(D,J)}(p,v)\ge \frac{c_{m=2t}}{2}\sqrt{\lambda_0(\tilde{D},\tilde{J},\tilde{u})}\,\frac{||v||}{\sqrt{|u(p)|}}=\frac{1}{2}\sqrt{\frac{2}{9k}}\,e^{-2t}\times\sqrt{\lambda_0(\tilde{D},\tilde{J},\tilde{u})}\,\frac{||v||}{\sqrt{|u(p)|}}.$$
Moreover, for all  $x\in{\tilde{D}}$ and $X\in\R^{2n}$,
\begin{eqnarray*}
\mathcal{L}^{\tilde{J}}_x\tilde{u}(X)=\mathcal{L}^J_{P_px}u(P_pX)&\ge& \lambda_0(D,J,u)\times||P_pX||^2\\ &\ge& \lambda_0(D,J,u)||X||^2/||P_p^{-1}||_0^2\ge\lambda_0(D,J,u)||X||^2/4.
\end{eqnarray*}
Hence
$\sqrt{\lambda_0(\tilde{D},\tilde{J},\tilde{u})}\ge\frac{1}{2}\sqrt{\lambda_0(D,J,u)}$,
and we can set
$c'=\frac{1}{4}\sqrt{\frac{2}{9k}}$ (where $k$ is given by (\ref{cstek})).
\end{dem}

Therefore, we obtain a lower bound for the Kobayashi metric in a neighborhood of any point $p$ in an almost complex manifold, with an estimate of the size of this neighborhood depending on the choice of the coordinate system:

\begin{cor}
Let $(M,J)$ be an almost complex manifold and $z:U\to \B^{2n}$ a local coordinate system such that $z_*J(0)$ is the standard complex structure. Assume that $D\subset M$ is a domain and $u:\overline{D\cap U}\to ]-\infty;0[$ is a strictly $J$-plurisubharmonic function.\\
Set $t=\mathrm{min}\left( 1,\,\frac{\varepsilon'}{||z_*J||_{\mathcal{C}^1(\bar{\mathcal{\bbb}}^{2n})}}\right)$ (where $\ve'$ is given by Proposition \ref{prop:minorationK'}) and $U_t=z^{-1}(t \B^{2n})$. Then, for every $p\in D\cap U_t$ and $v\in T_pM$, we have
$$K_{(D\cap U_t,J)}(p,v)\ge
c'\,e^{-2t}\sqrt{\lambda_0(z(D\cap U_t),z_*J,u\circ z^{-1})}\times\frac{||dz_p(v)||}{\sqrt{|u(p)|}}.$$
\end{cor}

\begin{dem}
For any $t\in ]0;1]$ and $z_t:U_t\ni p\mapsto z(p)/t\in\B^{2n}$, one gets $(z_t)_*J_x=(z_*J)_{tx}$. Hence
$||(z_t)_*J-J_{st}||_{\mathcal{C}^1(\bar{\bbb}^{2n})}\le\varepsilon'$.
This allows us to apply Proposition \ref{prop:minorationK'} to the domain $z_t(D\cap U_t)$ equipped with the almost complex structure $(z_t)_*J$, and to the function $u\circ z_t^{-1}$. Then, for every $\tilde{p}\in z_t(D\cap U_t)$ and $\tilde{v}\in\R^{2n}$,
$$K_{(z_t(D\cap U_t),(z_t)_*J)}(\tilde{p},\tilde{v})\ge c'\,e^{-2t}\sqrt{\lambda_0(z_t(D\cap U_t),(z_t)_*J,u\circ z_t^{-1})}\frac{||\tilde{v}||}{\sqrt{u(z_t^{-1}(\tilde{p}))}}.$$
But ${\mathcal{L}^{(z_t)_*J}_x(u\circ z_t^{-1})}(X)=\mathcal{L}^{z_*J}_{tx}(u\circ
z^{-1})(tX)$, thus
$$\lambda_0(z_t(D\cap U_t),{(z_t)}_*J,u\circ z_t^{-1})\ge t^2\lambda_0(z(D\cap U_t),z_*J,u\circ z^{-1}),$$
 which gives the desired inequality with $\tilde{p}=z_t(p)$ and $\tilde{v}={dz_t}_p(v)=\frac{1}{t}\,{dz}_p(v)$.
\end{dem}

\subsection{Localization principle}
The Kobayashi metric verifies $K_{(D\cap U,J)}(q,v)\ge K_{(D,J)}(q,v)$. Proposition 3 of \cite{GS} gives a sort of converse, which is very useful while working with charts on manifolds:
$$K_{(D,J)}(q,v)\ge sK_{(D\cap U,J)}(q,v).$$
We are going to compute explicitely the constant $s$. The proof will also provide some explicit estimate of the size of pseudo-holomorphic discs.

\begin{lm}\label{lm:lemmeprec}
Let $D$ be a domain in an almost complex manifold $(M,J)$. Pick $p\in\bar{D}$ and $z:U\to\B$ a chart in a neighborhood of $p$ verifying $z_*J(p)=J_{st}$ and $||z_*J-J_{st}||_{\mathcal{C}^1(\bar{D})}\le\ve_{m=1}$. Assume that there exists a $\mathcal{C}^2$ regular map $u:\bar{D}\to ]-\infty;0[$, strictly $J$-plurisubharmonic in $D$.\\ 
Then there is a neighborhood $V\Subset U$ of $p$ verifying that for all $q\in D\cap V$ and $v\in T_qM$:
$$K_{(D,J)}(q,v)\ge N||dz_q(v)||\,,\quad\mathrm{with}\ N=\frac{e^{-1}}{\sqrt{k}}\sqrt{\frac{c}{|u(q)|}}$$
where $c>0$ is such that $u-c||z||^2$ is strictly $J$-plurisubharmonic. Here $k$ is defined by (\ref{cstek}).
\end{lm}

\begin{dem}
Let $\theta$, $r$, $A$, $B$ be as in Lemma \ref{lm:epsilonm}. The map
$$x\overset{u}{\mapsto}\ln\left(\theta\left(\frac{1}{r^2}||z(x)-z(q)||^2\right)\right)+A||z(x)-z(q)||+B\frac{1}{r^2}||z(x)-z(q)||^2$$
is strictly $J$-plurisubharmonic in $U$ for $q=p$, and hence for all $q$ in some neighborhood $V\Subset U$ of $p$. Note that $V$ depends only on the choice of $z$ (and $\theta$). Pick $\lambda> 1/r^2$ and $\tau=\lambda B/c$.
For any $q\in V$, we define $\Psi_q$ by
$$\Psi_q(x)=\left\{\begin{array}{l}
\theta\left(\frac{1}{r^2}||z(x)-z(q)||^2\right)\,\mathrm{exp}(A||z(x)-z(q)||)\,\mathrm{exp}(\tau u(x))\quad\mathrm{if}\ x\in D\cap U,\\
\mathrm{exp}(A+\tau u(x))\quad\mathrm{if}\ x\in D\setminus U.
\end{array}\right.$$
As soon as $0<\ve<B(\lambda-1/r^2)$, the map $\mathrm{ln}\,(\Psi_q)-\ve||z-z(q)||^2$ is $J$-plurisubharmonic in $D\cap U$, and thus $\Psi_q$ is $J$-plurisubharmonic in $D\cap U$. But $\Psi_q$ coincides with $\mathrm{exp}\,(A+\tau u)$ out of $U$, hence is globally $J$-plurisubharmonic in $D$.

Let $h:\D\to D$ be a $J$-holomorphic disc such that $h(0)=q\in V$ and $dh_0(\partial/\partial x)=v/\alpha$, where $v\in T_qM$ and $\alpha>0$. For $\zeta$ sufficiently small, 
$$h(\zeta)=q+dh_0(\zeta)+O(|\zeta|^2).$$
Set $\zeta=\zeta_1+i\zeta_2$. The $J$-holomorphy equation $dh_0\circ J_{st}=J\circ dh_0$ gives
$$dh_0(\zeta)=\zeta_1dh_0(\partial/\partial x)+\zeta_2Jdh_0(\partial/\partial x).$$
Let us consider the map
$$\vp(\zeta)=\frac{\Psi_q(h(\zeta))}{|\zeta|^2},$$
which is subharmonic in $\D\setminus\{0\}$. Since
$$\vp(\zeta)=\frac{||z\circ h(\zeta)-z(q)||^2}{r^2|\zeta|^2}\,\mathrm{exp}\,(A||z\circ h(\zeta)-z(q)||)\,\mathrm{exp}\,(\tau u(h(\zeta)))$$
for small $\zeta$, and
$$||dh_0(\zeta)||\le|\zeta|(||I+J||\cdot||dh_0(\partial/\partial x)||),$$
we obtain that $\underset{\zeta\to 0}{\mathrm{limsup}}\,\vp(\zeta)$ is finite. Moreover, with $\zeta_2=0$, this gives:
$$\underset{\zeta\to 0}{\mathrm{limsup}}\,\vp(\zeta)\ge \frac{||dh_0(\partial/\partial x)||^2}{r^2}\mathrm{exp}\,(-\lambda B|u(q)|/c).$$
The maximum principle applyed to a subharmonic extension of $\vp$ to $\D$ implies that for all $q\in D\cap V$ and $v\in T_qM$:
$$K_{(D,J)}(q,v)\ge\frac{||dz_qv||}{r}\mathrm{exp}\,(-(A+\lambda B|u(q)|/c)/2).$$
With $A\to 1$, $B=k(\theta)$, $\lambda\to 1/r^2$, one obtains
$$K_{(D,J)}(q,v)\ge\frac{\mathrm{exp}\,(-(1+k|u(q)|/(cr^2))/2)}{r}\,||dz_qv||$$
which is minimal for $r=\sqrt{k|u(q)|/c}$.
\end{dem}

As a corollary, we get:

\begin{prop}\label{prop:tailledisques}
Assume that the conditions of Lemma \ref{lm:lemmeprec} hold. Then every $J$-holomorphic disc $h:\D\to D$ such that $h(0)\in V$ verifies $h(s\D)\subset V$, where 
$$s=1-\mathrm{exp}\,(-N\dist(z\circ h(0)),\partial\B).$$
\end{prop}

\begin{dem}
Set $q=h(0)$ and $s=1-\mathrm{exp}\,(-N\dist(z\circ h(0)),\partial\B)$. Assume by contradiction that there exists $\zeta\in\D$ such that $w=h(s\zeta)\notin V$. We define $G=\{x\in V/\ ||z(x)-z(q)||<\delta\}$ where $\delta=\dist(z(q),\partial\B)$. By hypothesis, we have for every path $\gamma:[0;1]\to M$ between $q$ and $h(s\zeta)$:
\begin{eqnarray*}
d_{(M,J)}^K(h(0),h(s\zeta))=\int_0^1K_{(M,J)}(\gamma(t),\gamma'(t))\mathrm{d}t
&\ge&\int_{\gamma^{-1}(G)}N\times||(z\circ\gamma)'(t)||\mathrm{d}t\\
&=&N\dist(z(q),\partial\B).
\end{eqnarray*}
But one obtains easily, considering the pseudo-holomorphic disc $g(\zeta)=\zeta_0+\frac{v}{|v|}(1-|\zeta_0|)\zeta$, that
$$d_{(M,J)}^K(h(0),h(\zeta))\le d^K_{(\Delta,J_{st})}(0,\zeta)\le -\ln(1-|\zeta|).$$
 This gives a contradiction in view of the definition of $s$.
\end{dem}

\begin{cor}
Define the neighborhood $V$ and the constant $s$ as previously. Then for all $q\in D\cap V$ and $v\in T_qM$:
$$K_{(D\cap U,J)}(q,v)\ge K_{(D,J)}(q,v)\ge sK_{(D\cap U,J)}(q,v).$$
\end{cor}

\section{$J$-holomorphic discs attached to a maximal totally real submanifold}
We already know (see \cite{CGS}) that if $J$ is smooth, then any analytic disc attached to a totally real sumanifold is also smooth up to the boundary. Here, we are looking for a more precise result: if the almost complex structure is only of class $\mathcal{C}^r$, what regularity can we expect? 

The first step is to get the $\frac{1}{2}$-Hölderian regularity up to the boundary.

\subsection{Estimate of {\boldmath $\frac{1}{2}$\unboldmath}-Hölderian norm}
\begin{lm}\label{lm:inegalitemoyenne}
Let $\phi$ be a non-negative subharmonic function on the unit disc $\Delta$, continuous up to the boundary, which vanishes on the upper half-circle.
Given $\alpha\in ]0;\frac{\pi}{2}[$, let us denote by $W_\alpha$ the angular sector $\{re^{i\theta}/\ 0<r\le 1,\,\alpha<\theta<\pi-\alpha\}$. Then
$$\forall\zeta\in
W_\alpha,\ 
\phi(\zeta)\le\left(\frac{1}{(\sin\alpha)^2}\times\int_\pi^{2\pi}\phi(e^{i\theta})\frac{d\theta}{\pi}\right)\times(1-|\zeta|).$$
\end{lm}

\begin{dem}
It is immediate if $|\zeta|=1$. Therefore, we assume that $\zeta=re^{it}\in W_\alpha$ for some $r<1$.
Since $\phi$ is subharmonic on $\Delta$ and continuous on $\overline{\D}$, we get 
$$\forall \zeta=re^{it}\in\Delta,\ \phi(\zeta)\le\frac{1}{2\pi}\int_0^{2\pi}\phi(e^{i\theta})\,\frac{1-r^2}{|e^{i\theta}-re^{it}|^2}\,d\theta.$$
Given $\alpha\in
]0;\frac{\pi}{2}[,\ \theta\in [\pi;2\pi]$ and $t\in
[\alpha;\pi-\alpha]$, we have $\alpha\le \theta-t\le\pi-\alpha$, and so
$r^2-2r\cos(\theta-t)+1\ge r^2-2r\cos\alpha+1\ge 1-(\cos\alpha)^2$. This gives the inequality. 
\end{dem}

As in the standard case, the $\frac{1}{2}$-Hölderian regularity and the estimate of the associated norm come from an estimate of the differential, depending on the square root of the distance to the boundary.

\begin{lm}
Let $D$ be a bounded domain in an almost complex manifold $(M,J)$, $\rho\in\mathcal{C}^2(\bar{D},\R)$ some strictly plurisubharmonic function and $h:\D\to D$ a $J$-holomorphic disc, continuous up to the boundary. We assume that $\rho\circ h\ge 0$ on $\Delta$ and vanishes on the upper half-circle $\gamma=\{e^{i\theta}/\,0<\theta<\pi\}$.

Pick $a\in\gamma$ and let $(z,U)$ be a local coordinate system such that $h(a)\in U$ and $||z_*J-J_{st}||_{\mathcal{C}^1(\bar{\B})}\le\varepsilon'$. Then, there exist a universal constant c'' and a neighborhood $V$ of $a$ in $\Delta\cap h^{-1}(U)$
such that
$$\forall\zeta\in
V,\ |||d(z\circ h)_\zeta|||\le c''\,\frac{1}{\mathrm{Im}\,a}\,\sqrt{\frac{\int_0^{2\pi}\rho\circ
  h(e^{i\theta})d\theta}{\lambda_0(z(D\cap
  U),z_*J,\rho\circ z^{-1})}}\times\frac{1}{\sqrt{1-|\zeta|}}\,.$$
\end{lm}

\begin{dem}
We use the notations of Lemma \ref{lm:inegalitemoyenne}. Let $\delta>0$ be such that 
$\Omega_\delta=\Delta\cap(a+\delta\Delta)$ is included in 
$W_\alpha$ for $\alpha=\mathrm{Arcsin}(\Im a/2)$. Reducing $\delta$ if necessary, we may assume by continuity of $h$ that 
$h(\Omega_\delta)\subset U$. Thus, Lemma \ref{lm:inegalitemoyenne} gives:
$$\forall\zeta\in\Omega_\delta,\ \rho\circ
h(\zeta)\le\kappa(1-|\zeta|)\quad \mathrm{where}\quad \kappa=\left(\frac{4}{(\Im a)^2}\times\int_0^{2\pi}\rho\circ
  h(e^{i\theta})d\theta\right).$$
Pick $\zeta_0\in\Omega_{\delta/2}$ and set $l=1-|\zeta_0|$: then $\zeta_0+l\Delta\subset\Omega_\delta$. For 
$D_l=\{q\in D/\ \rho(q)<2\kappa l\}$, we obtain $\rho\circ h(\zeta_0+l\Delta)\subset D_l$.\vspace*{0.1cm}\\
The function $u_l:w\mapsto\rho\circ z^{-1}(w)-2\kappa l$
is negative and strictly $z_*J$-plurisubharmonic in $z(D_l\cap U)\supset z\circ 
h(\zeta_0+l\Delta)$. Let us define
$$g_l:\Delta\ni\zeta\mapsto z\circ h(\zeta_0+l\zeta)\in z(D_l\cap
U).$$
Since the disc $g_l$ is $z_*J$-holomorphic, we get, by the decreasing property of the Kobayashi metric:
$$
|\tau|\ge K_{(z(D_l\cap U),z_*J)}(z\circ h(\zeta_0),l\times d(z\circ
h)_{\zeta_0}(\tau)).
$$
By Proposition \ref{prop:minorationK'} applyed to $z(D_l\cap U)$, $z_*J$ and $u_l$, we obtain
$$\forall p\in z(D_l\cap U),\forall v\in \R^{2n},\ K_{(z(D_l\cap U),z_*J)}(p,v)\ge
c'e^{-2}\,\sqrt{\lambda_0(z(D_l\cap
  U),z_*J,u_l)}\times\frac{||v||}{\sqrt{|u_l(p)|}}.$$
With $p=z\circ h(\zeta_0)$ and $v=l\times
d(z\circ h)_{\zeta_0}(\tau)$, this gives
$$|\tau|\ge c'e^{-2}\,\sqrt{\lambda_0(z(D_l\cap
  U),z_*J,u_l)}\times\frac{||l\times d(z\circ
  h)_{\zeta_0}(\tau)||}{\sqrt{|\rho\circ h(\zeta_0)-2\kappa l|}}\,.$$
Hence
$$|||d(z\circ h)_{\zeta_0}|||\le\frac{\sqrt{2\kappa l}}{ c'e^{-2}\,\sqrt{\lambda_0(z(D_l\cap
  U),z_*J,u_l)}\times l}$$
and one can choose $c''=\frac{2e^2\sqrt{2}}{c'}$.
\end{dem}

One obviously obtains a similar result by replacing $\D$ by the upper half-disc $\D^+=\{\zeta\in\D/\ \Im\zeta>0\}$ and the upper half-circle by $]-1;1[$. Then, the estimate of the differential and Hardy-Littlewood give the local $\frac{1}{2}$-Hölderian continuity of $z\circ h$ on $\Delta^+\cup ]-1;1[$:

\begin{cor}\label{cor:csqKob}
Let $D$ be a bounded domain in an almost complex manifold $(M,J)$, $\rho\in\mathcal{C}^2(\bar{D},\R)$ some strictly plurisubharmonic function and $h:\D^+\to D$ a $J$-holomorphic map. We assume that $h$ extends continuously on $\D^+\cup]-1;1[$, such that 
$$\rho\circ h\ge 0\ \mathrm{on}\ \Delta^+\ \mathrm{and}\ \rho\circ h_{|]-1;1[}\equiv 0.$$
Pick $a\in ]-1;1[$ and let $(z,U)$ be a local coordinate system in a neighborhood of $h(a)$ such that $||z_*J-J_{st}||_{\mathcal{C}^1(\bar{\B})}\le\varepsilon'$. Then, there exist a universal constant c and a neighborhood $W$ of $a$ in $\Delta\cap h^{-1}(U)$
such that
$$\forall\zeta,\zeta'\in
W,\ ||z\circ h(\zeta)-z\circ h(\zeta')||\le c\,\frac{1}{1-|a|}\,\sqrt{\frac{\int_0^{2\pi}\rho\circ
  h(e^{i\theta})d\theta}{\lambda_0(z(D\cap
  U),z_*J,\rho\circ z^{-1})}}\times|\zeta-\zeta'|^{1/2}.$$
\end{cor}

We want to apply this result to the case of a pseudo-holomorphic disc attached to a maximal totally real submanifold. Let us recall that any maximal totally real submanifold in an almost complex manifold $(M^{2n},J)$ can be described as the zero set of some non-negative function $\rho$ as regular as $E$, strictly $J$-plurisubharmonic.
Indeed, if $E$ is defined on $U$ by $r_1=\hdots=r_n=0$ such that $dr_1\wedge\hdots\wedge dr_n$ does not vanish, $\rho=\sum_{i=1}^nr_i^2$ is strictly $J$-plurisubharmonic on $U$. 
For $(U^\alpha,\phi_\alpha)$ being a partition of unity, denote by $\rho_\alpha$ the function constructed as previously on the open set $U^\alpha$, and set $\rho=\sum_\alpha\rho_\alpha\phi_\alpha$. Then
 $\nabla\rho_{|E}\equiv 0$ and $$D_{|E}=\left(\left.\frac{\partial^2\rho}{\partial x_k\partial x_l}\right|_{E}\right)_{k,l}=2\sum_\alpha\phi_\alpha \left(\left.\frac{\partial^2\rho_\alpha}{\partial x_k\partial x_l}\right|_{E}\right)_{k,l}$$
is positive semi-definite. Moreover, $\tsp XDX=2\sum_{i=1}^n||\tsp \nabla r_{\alpha,i}\cdot X||^2=0$ on $E$ if and only if for all $\alpha$, $i$, $dr_{\alpha,i}(X)=0$: that is, $X\in TE$. 
It shows therefore that for any $p\in E$,
$$\mathcal{L}^J_{p}\rho(X)=\tsp XD(p)X+\tsp (J(p)X)D(p)(J(p)X)\ge 0$$
and 
$$\mathcal{L}^J_{p}\rho(X)=0\Longleftrightarrow \left(X\in TE\ \mathrm{et}\ JX\in TE\right)\Longleftrightarrow X=0$$
since $E$ is totally real.
The Levi form of $\rho$ is thus positive definite on $E$, hence in a neighborhood of $E$.

Moreover, we get that, if $\pi_{y_i}$ denotes the canonical projection from $\R^{2n}$ on the $y_i$-axis, the function $\rho=\sum_{i=1}^n(\pi_{y_i}\circ z)^2$ extends in a non-negative and strictly $J$-plurisubharmonic function, such that $E=\{\rho=0\}$. Thus, we have: $$\forall\zeta\in\D^+,\ |\rho\circ h(\zeta)|\le||z\circ h(\zeta)||^2\le||z\circ h||_\infty^2.$$
Hence, Corollary \ref{cor:csqKob} gives a neighborhood $W$ of $a$ such that 
$$\forall\zeta,\zeta'\in W,\ ||z\circ h(\zeta)-z\circ h(\zeta')||\le\frac{c}{1-|a|}\,\times\,\sqrt{\frac{\int_0^\pi\rho\circ
    h(e^{i\theta})d\theta}{\lambda_0(z(U\cap V),z_*J,\rho\circ z^{-1})}}\,\times|\zeta-\zeta'|^{1/2}.$$

We introduce the following notation:

\begin{defi}
The \textit{minimal $J$-curvature} of $E$ is
$$\lambda^J_E=\underset{i}{\mathrm{min}}\,\lambda_0(\B^{2n},{z_i}_*J,u),$$ 
where $u(x_1,\hdots,x_n,y_1,\hdots,y_n)=y_1^2+\hdots+y_n^2$ and $\lambda_0(\B^{2n},{z_i}_*J,u)$ is the smallest eigenvalue of the Levi form of $u$. 
\end{defi}

We have proved:

\begin{prop}\label{prop:normeC1/2}
Let $(M,J)$ be an almost complex manifold, $E$ a maximal totally real submanifold, $p\in E$ and $(z,U)$ a local coordinate system supposed to be $(\ve',E)$-tamed in a neighborhood of $p$.
Assume that $h$ is continuous from $\Delta^+\cup ]-1;1[$ to $M$, $J$-holomorphic in $\Delta^+$, and verifies $h(]-1;1[)\subset
  E$.\\
Then for all $a\in ]-1;1[\cap h^{-1}(U)$, there exists a neighborhood $W\ni a$ in $U$ such that $h$ is Hölderian $\frac{1}{2}$-continuous on $W$, and
$$\forall \zeta,\zeta'\in W,\ ||z\circ h(\zeta)-z\circ
h(\zeta')||\le\frac{\tilde{c}}{1-|a|}\,\times\,\frac{||z\circ h||_\infty}{\sqrt{\lambda_E^J}}\,\times|\zeta-\zeta'|^{1/2}$$
where $\tilde{c}$ is a universal constant, and $\lambda_E^J$ is the minimal $J$-curvature of $E$.
\end{prop}

\subsection{Estimates of Hölderian norms}
This estimate of $\frac{1}{2}$-Hölderian norm automatically leads to an estimate of Hölderian norms of upper degree. We will use the following proposition:

\begin{prop}\label{Siko}
{\bf (see \cite{S})} Pick $\alpha\in ]0;1[$, and let $\Omega$ be a domain in $\Delta$ and $K$ some relatively compact subset in $\Omega$.
There exist $\delta_\alpha>0$ and $\Lambda(\alpha,K)>0$ such that for every 
$q:\mathcal{\B}^n\to\mathrm{End}_{\rrr}(\C^n)$ of class $\mathcal{C}^{\alpha}$ verifying $||q||_{\alpha}\le\delta_\alpha$, every differentiable function $h:\Omega\to\mathcal{\B}^n$ verifying
$\bar{\partial}h+q\circ h\,\times\partial h=0$ is of class
$\mathcal{C}^{1+\alpha}$ in $K$ and 
$$||h_{|K}||_{1+\alpha}\le \Lambda(\alpha,K)\times ||h_{|K}||_{1/2}.$$
\end{prop}

\begin{dem}
We follow the proof of \cite{S}, choosing the function $\rho_1$
strictly supported in the interior of $\Omega$. Let us denote by $T$ and $P$ the operators defined by 
$$Pg(z)=\frac{1}{2\pi i}\int\int_D\frac{g(\zeta)}{\zeta-z}d\zeta\wedge
d\bar{\zeta}\quad\mathrm{and}\quad Tg(z)=p.v.\left(\in\int_D\frac{g(\zeta)}{(\zeta-z)^2}d\zeta\wedge d\bar{\zeta}\right).$$
We fix $K'$ and $K''$ two compact subsets such that $K\Subset
K''\Subset K'\Subset D$, and $\rho_1$, $\rho_2$, $\rho_3$ smooth functions from $\Delta$ to $[0;1]$ verifying
$$\mathrm{Supp}\,(\rho_1)\subset \mathrm{int}\,\Omega\ ,\ \mathrm{Supp}\,(\rho_2)\subset K'\ ,\ \mathrm{Supp}\,(\rho_3)\subset K''$$
and 
$${\rho_1}_{|K'}\equiv 1\ ,\ {\rho_2}_{|K''}\equiv 1\ ,\ {\rho_3}_{|K}\equiv 1.$$
Set
$M(K,D,\alpha)=max(||\rho_1||_{\mathcal{C}^{1+\alpha}},||\rho_2||_{\mathcal{C}^{1+\alpha}},||\rho_3||_{\mathcal{C}^{1+\alpha}})$.
Then, the constants $\Lambda$ and $\delta$ depend only on $|||T|||_{L^3(D)}$,
$|||P|||_{L^3(D)\to\mathcal{C}^{1/3}(D)}$,
$|||T|||_{\mathcal{C}^{\alpha/3}}$ and $M(K,D,\alpha)$.
\end{dem}

Using geometric bootstrap's method, we will obtain simultaneously the regularity and the estimates. The idea (see \cite{GSmodele}) is to symmetrize discs along the piece of their boundary which is attached to $E$. We then obtain a new disc, pseudo-holomorphic on all the domain. Since the pseudo-holomorphy condition can be expressed as an elliptic equation, we get the result by induction.

We assume that the hypotheses of Theorem \ref{p4.7-CGS} are satisfied. Let us define $D^+=h^{-1}(U)$, $D^-=\overline{h^{-1}(U)}$, $\delta=\overline{D^+}\cap\overline{D^-}$ and
$$g(\zeta)=\begin{array}{|l}
{z}\circ h(\zeta)\ \mathrm{if}\ \zeta\in D^+\\
\overline{{z}\circ h(\bar{\zeta})}\ \mathrm{if}\ \zeta\in D^-
\end{array}\,,$$
where $z$ is a chart as in Lemma \ref{lm:chgtcarte}.
The map $g$ is continuous on $D=D^+\cup
D^-\cup\delta$. Moreover, $g$ verifies in
$D$ the equation: 
$$\bar{\partial}g+A(\cdot)\partial g=0,$$
where $A(\zeta)=Q(g(\zeta))$ if $\zeta\in D^+\cup\delta$, and $A(\zeta)=\overline{Q(g(\bar{\zeta}))}$ if $\zeta\in D^-$, and $Q$ is given by Lemma \ref{lm:chgtcarte}. Thus, setting $\mathrm{min}(r,r')=k+\alpha$ with $k\in\N^*$ and $0<\alpha<1$, $Q$ is of class $\mathcal{C}^{k-1+\alpha}$.\\ 

Let $K$ be a compact subset in $(D^+\cup\delta)$.
If $k=1$, $A$ is of class $\mathcal{C}^{\alpha/2}$ since $g$ is Hölder continuous with exponent $1/2$. By Proposition \ref{Siko}, we thus get that $g$ is in fact of class $\mathcal{C}^{1+\alpha/2}$. But this implies that $A$ is of class $\mathcal{C}^\alpha$: hence by the same argument $g$ is of class $\mathcal{C}^{1+\alpha}$ in $K$ and 
\begin{equation}\label{k=1}
||z\circ h||_{\mathcal{C}^{1+\alpha}(K)}=||g||_{\mathcal{C}^{1+\alpha}(K)}\le \Lambda(\alpha,K)||g||_{\mathcal{C}^{1/2}(K)}=\Lambda(\alpha,K)||z\circ h||_{\mathcal{C}^{1/2}(K)}.
\end{equation}

Assume $k\ge 2$. Then $A$ is of class $\mathcal{C}^{1/2}$ and thus $g$ is in fact of class $\mathcal{C}^{1+1/2}$. Hence $A$ is $\mathcal{C}^1$ regular and $g$ is of class $\mathcal{C}^{2-0}$ in $K$.
 The \textsl{a priori} estimates of Proposition \ref{Siko} (see also \cite{Vekua}) apply to $g$: for all $0<\beta<1$,
$$||z\circ h||_{\mathcal{C}^{1+\beta}(K)}=||g||_{\mathcal{C}^{1+\beta}(K)}\le \Lambda(\beta,K)||g||_{\mathcal{C}^{1/2}(K)}=\Lambda(\beta,K)||z\circ h||_{\mathcal{C}^{1/2}(K)}.$$
According to Lemma \ref{lm:relevement}, the previous argument applies to $H^1=(h,dh)$. If $k-1\ge 2$, then $H^1$ is of class $\mathcal{C}^{1+\beta}$ in $K$ for all $0<\beta<1$. In particular, 
$||z^c\circ H^1||_{\mathcal{C}^{1+1/2}(K)}\le \Lambda(1/2,K)||z^c\circ
H^1||_{\mathcal{C}^{1/2}(K)},$
that is,
$$||z\circ h||_{\mathcal{C}^{2+1/2}(K)}\le \Lambda(1/2,K)||g||_{\mathcal{C}^{1+1/2}(K)}\le\Lambda(1/2,K)^2||{z}\circ
h||_{\mathcal{C}^{1/2}(K)}.$$
We obtain by iteration that $H^{k-2}=(g,dg,\hdots,d^{k-2}g)$ is of class $\mathcal{C}^{2-0}$ in $K$. Hence $z\circ h$ is of class $\mathcal{C}^{k-1+1/2}$ in $K$, and
$$||z\circ h||_{\mathcal{C}^{k-1+1/2}(K)}\le C(k,K)||z\circ h||_{\mathcal{C}^{1/2}(K)}.$$
Following the previous argument (case $k=1$), we get that $H^{k-1}$ is of class $\mathcal{C}^{1+\alpha}$. Consequently, $z\circ h$ is of class $\mathcal{C}^{k+\alpha}$ and in view of (\ref{k=1}), 
$$||z\circ h||_{\mathcal{C}^{k+\alpha}(K)}\le C(K,k)\Lambda(\alpha,K)||z\circ h||_{\mathcal{C}^{1/2}(K)}.$$
Moreover,
 $$||{z}\circ h||_{\mathcal{C}^{1/2}(K)}\le||{z}\circ h||_\infty\left(1+\frac{c(K)}{\sqrt{\lambda_E^J}}\right)$$
by means of Proposition \ref{prop:normeC1/2}. This exactly means that $h$ is of class $\mathcal{C}^{\mathrm{min}(r,r')}$ and that 
$$||h||_{\mathcal{C}^{\mathrm{min}(r,r')}(K)}\le c(r,K)||h||_\infty\left(1+\frac{c(K)}{\sqrt{\lambda^J_E}}\right)\,.$$
Hence we have proved the second part of Theorem \ref{p4.7-CGS}.
The construction of an atlas $(\ve',E)$-tamed also shows that the choice of another atlas does not modify the estimate.

\subsection{Regularity on the edge of the wedge}
As in \cite{CGS}, the regularity and the estimates given by Theorem \ref{p4.7-CGS} for pseudo-holomorphic discs provide the regularity and similar estimates for pseudo-holomorphic maps defined on a wedge. We begin with some preliminary remarks.

Let $\Omega\subset M^{2n}$ be a domain. Suppose that $N\subset\Omega$ is defined by $r_1=\hdots=r_n=0$, where
$dr_1\wedge\hdots\wedge dr_n$ does not vanish on $\Omega$. Let us denote by $W(\Omega,N)=\{z\in\Omega/\ \forall
1\le j\le n,\ r_j(z)<0\}$ the wedge of edge $N$.\\
The restricted wedge 
$W_\delta(\Omega,N)=\{z\in\Omega/\ \forall 1\le j\le n,\
r_j(z)-\delta\sum_{k\not=j}r_k(z)<0\}$ is included in $W(\Omega,N)$ for all $0<\delta<\frac{1}{n-1}$. Note that $W_\delta(\Omega,N)$ and $W(\Omega,N)$ 
are non-empty open subsets in $\Omega$, and that $N$ is included in their boundary. 

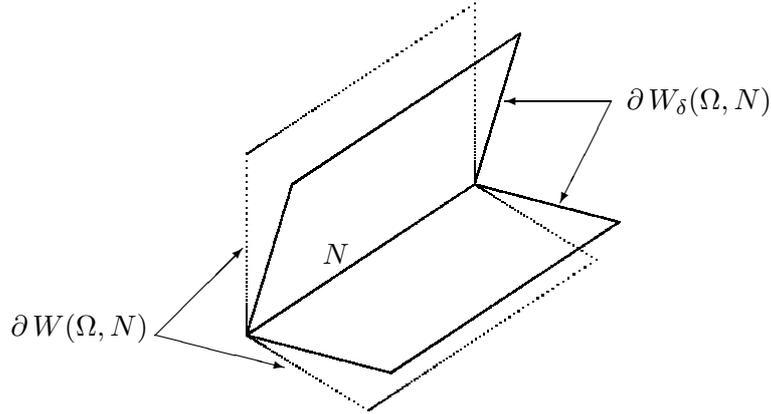
\begin{picture}(15,6.5)
\qbezier(4.5,1.8)(4.5,1.8)(7.5,3.8)
\qbezier(4.5,1.8)(4.5,1.8)(6.4,1.3)
\qbezier[50](4.5,1.8)(4.5,1.8)(6.1,0.8)
\qbezier(6.4,1.3)(6.4,1.3)(9.4,3.3)
\qbezier(7.5,3.8)(7.5,3.8)(9.4,3.3)
\qbezier(4.5,1.8)(4.5,1.8)(5.1,3.8)
\qbezier(5.1,3.8)(5.1,3.8)(8.1,5.8)
\qbezier(7.5,3.8)(7.5,3.8)(8.1,5.8)
\qbezier[50](4.5,1.8)(4.5,1.8)(4.5,4.2)
\qbezier[80](4.5,4.2)(4.5,4.2)(7.5,6.2)
\qbezier[50](7.5,3.8)(7.5,3.8)(7.5,6.2)
\qbezier[80](6.1,0.8)(6.1,0.8)(9.1,2.8)
\qbezier[50](9.1,2.8)(7.5,3.8)(7.5,3.8)

\put(5.5,2.75){$N$}
\put(1.4,1.8){$\partial\, W(\Omega,N)$}
\put(3.3,1.8){\vector(1,1){1.15}}
\put(3.3,1.8){\vector(4,-1){1.7}}
\put(9.5,4.8){$\partial\, W_\delta(\Omega,N)$}
\put(9.3,4.9){\vector(-1,0){1.4}}
\put(9.3,4.9){\vector(-1,-2){0.67}}
\end{picture}

Let us prove Corollary \ref{cor:wedge}.
Assume that the conditions of the Corollary hold. In view of Lemma 5.2 in \cite{CGS}, there exists a parametrized family $(h_t)_{t\in\rrr^{2n}}$ of $J$-holomorphic half-discs, smoothly depending on $t$,
such that $W_\delta(\Omega,N)\subset\bigcup_th_t(\Delta^+)$ and
$$\forall t,\ h_t(]-1;1[)\subset N\ \ \mathrm{and}\ \ h_t(\Delta^+)\subset W(\Omega,N).
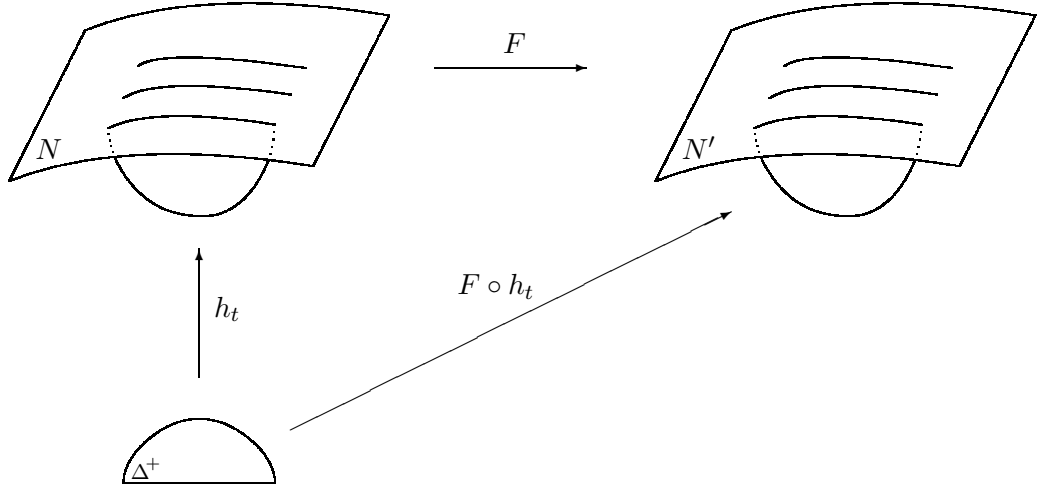
$$

\begin{picture}(14,8)
\qbezier(0.5,4.7)(2,5.3)(4.5,4.9)
\qbezier(1.5,6.7)(3,7.3)(5.5,6.9)
\qbezier(0.5,4.7)(0.5,4.7)(1.5,6.7)
\qbezier(4.5,4.9)(4.5,4.9)(5.5,6.9)
\put(0.85,5){$N$}

\qbezier(9,4.7)(10.5,5.3)(13,4.9)
\qbezier(10,6.7)(11.5,7.3)(14,6.9)
\qbezier(9,4.7)(9,4.7)(10,6.7)
\qbezier(13,4.9)(13,4.9)(14,6.9)
\put(9.35,5){$N'$}

\qbezier(1.8,5.4)(2.4,5.7)(4,5.45)
\qbezier[5](1.8,5.4)(1.82,5.2)(1.89,5.01)
\qbezier(1.89,5.01)(2.3,4.15)(3.2,4.25)
\qbezier(3.2,4.25)(3.65,4.35)(3.91,4.97)
\qbezier[5](3.91,4.97)(3.98,5.2)(4,5.45)

\qbezier(2,5.8)(2.45,6.1)(4.2,5.85)
\qbezier(2.2,6.23)(2.55,6.47)(4.4,6.2)
\qbezier(1.8,5.4)(2.4,5.7)(4,5.45)

\qbezier(10.3,5.4)(10.9,5.7)(12.5,5.45)
\qbezier[5](10.3,5.4)(10.32,5.2)(10.39,5.01)
\qbezier(10.39,5.01)(10.8,4.15)(11.7,4.25)
\qbezier(11.7,4.25)(12.15,4.35)(12.41,4.97)
\qbezier[5](12.41,4.97)(12.48,5.2)(12.5,5.45)

\qbezier(10.5,5.8)(10.95,6.1)(12.7,5.85)
\qbezier(10.7,6.23)(11.05,6.47)(12.9,6.2)
\qbezier(10.3,5.4)(10.9,5.7)(12.5,5.45)

\qbezier(2,0.7)(2.03,1.1)(2.5,1.4)
\qbezier(2.5,1.4)(3,1.7)(3.5,1.4)
\qbezier(3.5,1.4)(3.97,1.1)(4,0.7)
\put(2,0.7){\line(1,0){2}}
\put(2.1,0.9){$_{\D{\!\!}^{\,^+}}$}

\put(3,2.1){\vector(0,1){1.7}}
\put(3.2,2.9){$h_t$}

\put(6.1,6.2){\vector(1,0){2}}
\put(7,6.4){$F$}

\put(4.2,1.4){\vector(2,1){5.8}}
\put(6.4,3.2){$F\circ h_t$}
\end{picture}

Theorem \ref{p4.7-CGS}
shows that the functions $F\circ h_t$ are locally of class $\mathcal{C}^{r'}$ in  $\Delta^+\cup
  ]-1;1[$. It also gives that for every compact set $K\subset\D^+\cup]-1;1[$: 
$$||F\circ h_t||_{\mathcal{C}^{r'}(K)}\le c(r',K)||F\circ h_t||_\infty\left(1+\frac{c(K)}{\sqrt{\lambda^{J'}_{N'}}}\right).$$
Hence, the $\mathcal{C}^{r'}$-Hölderian norms of the functions $F\circ h_t$ are uniformly bounded.
%
%
 Since the curves $h_t(]-1;1[)$ form
a family of transversal $\mathcal{C}^{r}$-foliations of $N$, one gets by the separated regularity principle \cite{Tregsep} that $F$ is of class $\mathcal{C}^{s}$. Moreover, 
the uniform upper bound for the $\mathcal{C}^{s}$-Hölderian norms of the $F\circ h_t$ also bounds the $\mathcal{C}^{s}$-Hölderian norm of $F$ up to a multiplicative universal constant \cite{Meyer}. This concludes the proof.

\appendix
\section{Small perturbations of the standard structure}
\subsection{Vector fields and $J$-holomorphic forms}
Let $J=\left(\begin{array}{cc}A&B\\C&D\end{array}\right)$ be a small perturbation of the standard structure as in Section \ref{s:SmallPerturbation}. For all $1\le j\le n$, we set
$$\frac{\partial^J}{\partial z_j}=\frac{1}{2}(I_{2n}-iJ)\frac{\partial}{\partial x_j}\quad\mathrm{and}\quad\frac{\partial^J}{\partial \bar{z}_j}=\frac{1}{2}(I_{2n}+iJ)\frac{\partial}{\partial x_j}.$$
Then $(\frac{\partial^J}{\partial z_1},\hdots,\ \frac{\partial^J}{\partial z_n},\,\frac{\partial^J}{\partial \bar{z}_1},\hdots,\ \frac{\partial^J}{\partial \bar{z}_n})$ is a $\C$-basis of $TM\otimes\,\C$, and consequently its dual basis $(d^Jz_1,\hdots,\,d^Jz_n,\,d^J\bar{z}_1,\hdots,\,d^J\bar{z}_n)$ is given in local coordinates by the matrix
$$\left(\begin{array}{cc}
I_n&I_n\\
i\tsp((I_n+iA)C^{-1})&-i\tsp((I_n-iA)C^{-1})
\end{array}\right).$$
Hence $d^Jz=dx+Pdy=dz+Hdy$ and $d^J\bar{z}=dx+\overline{P}dy=d\bar{z}+\overline{H}dy$,
where 
$$P=i\tsp\left((I_n+iA)C^{-1}\right)=iI_n+H$$ 
with ``small'' $H=R+iS$. Thus 
$C=(S+I_{n})^{-1}$ and $A=-R(S+I_n)^{-1}$, which defines $J$ uniquely by Section \ref{s:SmallPerturbation}:

\begin{prop}\label{prop:unicite}
For every family of differential forms $\omega=dz+Hdy$, where $H$ is near $0_n$, there is a unique almost complex structure $J$ near $J_{st}$ such that $\omega$ is a basis of $(1,0)$-forms with respect to $J$.
Conversely, every small perturbation of the standard structure admits a basis of $(1,0)$-forms of this type.
\end{prop}     

Let $h:\D\to\R^{2n}$ be a $\mathcal{C}^1$ regular map. We denote by $(u,v)$ the canonical coordinates in $\R^2$. Then 
 $h$ is $J$-holomorphic if and only if 
 $\frac{\partial h}{\partial u}+i\frac{\partial h}{\partial v}$ is a $(0,1)$-vector field, where 
 $i$ is the standard complex structure on $TM\otimes\C$. Assume that we are in the previous situation:

$$h\ \mathrm{is}\ J-\mathrm{holomorphic}\ \Longleftrightarrow (I_n\quad\tsp P)\left(\frac{\partial h}{\partial u}+i\frac{\partial h}{\partial v}\right)=0.$$
By setting
$$\frac{\partial h}{\partial \zeta}=\frac{1}{2}\left(\frac{\partial h}{\partial u}-J_{st}\frac{\partial h}{\partial v}\right)\quad\mathrm{and}\quad\frac{\partial h}{\partial \bar{\zeta}}=\frac{1}{2}\left(\frac{\partial h}{\partial u}+J_{st}\frac{\partial h}{\partial v}\right),$$
and separating real and imaginary parts, this gives that 
 $h$ is a $J$-holomorphic disc if and only if 
\begin{equation}\label{eqndholomorphie}
\forall\zeta\in\Delta,\ \left(\begin{array}{cc}
I-C^{-1}&-AC^{-1}\\
-AC^{-1}&-I+C^{-1}
\end{array}\right)_{|h(\zeta)}\,\frac{\partial h}{\partial {\zeta}}+
\left(\begin{array}{cc}
I+C^{-1}&-AC^{-1}\\
AC^{-1}&I+C^{-1}
\end{array}\right)_{|h(\zeta)}\,\frac{\partial h}{\partial \bar{\zeta}}=0\,.
\end{equation}

\subsection{Proof of Lemma \ref{lm:chgtcarte}}
Let $(z',U')$ be some chart in a neighborhood of $p$ such that  
$z'(p)=0$, ${(z'_*J)}_0=J_{st}$. We also assume that $z'(U\cap E)=\R^n\cap \B^{2n}$ when $p\in E$. For any $t>0$, we denote by $d_t:q\mapsto q/t$ the dilation defined on $\R^{2n}$, and by $z_t$ the map $z_t=d_t\circ z'$. For $t$ sufficiently small, $z=z_t$ and $U=z_t^{-1}(\B)$ verify $||z_*J-J_{st}||_{\mathcal{C}^1}\le\varepsilon$. Thus we suppose that $M=\B^{2n}$, $E=\R^n\cap\B^{2n}$ and in view of (\ref{petiteperturbation})  
$$J=\left(\begin{array}{cc}A&-(I_n+A^2)C^{-1}\\C&-CAC^{-1}\end{array}\right)\quad\mathrm{where}\quad ||A||_{\mathcal{C}^1}\le\ve''\quad\mathrm{and}\quad ||C-I_n||_{\mathcal{C}^1}\le\ve''.$$
Set 
$$\phi(x,y)=\left(\begin{array}{cc} I_n&-A(x,y)C(x,y)^{-1}\\ 0& C(x,y)^{-1}\end{array}\right)\left(\begin{array}{c}x\\y\end{array}\right)=\left(\begin{array}{l}x^*=x-AC^{-1}y\\ y^*=C^{-1}y\end{array}\right).$$
Since $$d\phi_{(x,y)}\leftrightarrow \left(\begin{array}{cc}I_n-\sum_{k=1}^n\left(\frac{\partial (AC^{-1})_{i,k}}{\partial x_j}\right)_{(i,j)}\,y_k\ &\ -\sum_{k=1}^n\left(\frac{\partial (AC^{-1})_{i,k}}{\partial y_j}\right)_{(i,j)}\,y_k-AC^{-1}\\
\sum_{k=1}^n\left(\frac{\partial (C^{-1})_{i,k}}{\partial x_j}\right)_{(i,j)}\,y_k\ &\ \sum_{k=1}^n\left(\frac{\partial (C^{-1})_{i,k}}{\partial y_j}\right)_{(i,j)}\,y_k+C^{-1}
\end{array}\right)$$
is near $I_{2n}$, there is a neighborhood $\tilde{U}$ of $p$ such that $\phi$ induces a local diffeomorphism of class $\mathcal{C}^r$ from $\tilde{U}$ to $\phi(\tilde{U})$. We set $\tilde{z}=\phi\circ z$: it remains to verify that the three conditions of the Lemma hold. 

Condition 1. is immediate. Moreover, 
 $d\phi_{(x,y)}\leftrightarrow \left(\begin{array}{cc}I_n& -A(x,y)C^{-1}(x,y)\\0&C^{-1}(x,y)\end{array}\right)+H(x,y)$ where
$$H(x,y)=\left(\begin{array}{cc}-\sum_{k=1}^n\left(\frac{\partial (AC^{-1})_{i,k}}{\partial x_j}\right)_{(i,j)}\,y_k & -\sum_{k=1}^n\left(\frac{\partial (AC^{-1})_{i,k}}{\partial y_j}\right)_{(i,j)}\,y_k\\
\sum_{k=1}^n\left(\frac{\partial (C^{-1})_{i,k}}{\partial x_j}\right)_{(i,j)}\,y_k&\sum_{k=1}^n\left(\frac{\partial (C^{-1})_{i,k}}{\partial y_j}\right)_{(i,j)}\,y_k
\end{array}\right).$$
Thus
$$d\phi^{-1}_{\phi(x,y)}=(d\phi_{(x,y)})^{-1}\leftrightarrow \left(I_{2n}+\left(\begin{array}{cc}I_n&A\\0&C\end{array}\right)\times H\right)^{-1}\times\left(\begin{array}{cc}I_n&A\\0&C\end{array}\right)$$
and
\begin{eqnarray*}
{(\phi_*J)_{(x^*,y^*)=\!\phi(x,y)}}&=&d\phi_{(x,y)}\circ J_{(x,y)}\circ d\phi^{-1}_{\phi(x,y)}\\
  &\leftrightarrow&\!\left(\begin{array}{cc}I_n&-AC^{-1}\\ 0&C^{-1}\end{array}\right)\!\left(\begin{array}{cc}A&-(I_n+A^2)C^{-1}\\ C&-CAC^{-1}\end{array}\right)\!\left(\begin{array}{cc}I_n&A\\ 0&C\end{array}\right)+\tilde{H}\\
 &=&\!J_{st}+\tilde{H}
\end{eqnarray*}
where $||H(x^*,y^*)||\le c\,\varepsilon''||y^*||$ for some constant $c$. This is Condition 2.

Finally, we obtain by (\ref{eqndholomorphie}) that the $\tilde{z}_*J$-holomorphy equation satisfied by discs $g=\tilde{z}\circ h$ is under the form:
$$\frac{\partial g}{\partial \bar{\zeta}}+Q\circ g(\zeta)\times\frac{\partial g}{\partial \zeta}=0,$$
where $Q(x^*,y^*)=(I_{2n}+O(||y^*||))^{-1}\times O(||y^*||)\le c\ ||y^*||$ for some constant $c$. Note that $Q$ has the regularity of $\tilde{z}_*J$.

\subsection{Proof of Lemma \ref{lm:relevement}}
Denote by $(x_1,\hdots,x_n,y_1,\hdots,y_n)$ the coordinates on $M$ given by $z$, and by $(x_1,\hdots,x_n,X_1,\hdots,X_n,y_1,\hdots,y_n,Y_1,\hdots,Y_n)$ the coordinates on $M^c$ given by $z^c$ (where $X_i,\ Y_i$ are fiber's coordinates). If  
$J=\left(\begin{array}{cc}A&B\\ C&D\end{array}\right)$, we set:
$$\alpha=\sum_{k=1}^n(X_k\frac{\partial A}{\partial x_k}+Y_k\frac{\partial A}{\partial y_k})\ ,\ \beta=\sum_{k=1}^n(X_k\frac{\partial B}{\partial x_k}+Y_k\frac{\partial B}{\partial y_k}),$$
$$\gamma=\sum_{k=1}^n(X_k\frac{\partial C}{\partial x_k}+Y_k\frac{\partial C}{\partial y_k})\ ,\ \delta=\sum_{k=1}^n(X_k\frac{\partial D}{\partial x_k}+Y_k\frac{\partial D}{\partial y_k}).$$
Then
$$J^c=\left(\begin{array}{cccc}
A&0&B&0\\
\alpha&A&\beta &B\\
C&0&D&0\\
\gamma&C&\delta&D
\end{array}\right)$$
defines an almost complex structure on $M^c$ for which $E^c=TE$ is totally real \cite{GSmodele}. The change of variables given by Lemma \ref{lm:chgtcarte} is  
$$\phi(x,y)=\left(\begin{array}{cc}I_n&-AC^{-1}\\0&C^{-1}\end{array}\right)\left(\begin{array}{c}x\\ y\end{array}\right)$$ 
and the one for $M^c$ is
$$\phi^c(x,X,y,Y)=\left(\begin{array}{cc}
I_{2n}&-\left(\begin{array}{cc}A&0\\ \alpha&A\end{array}\right)\left(\begin{array}{cc}C&0\\ \gamma&C\end{array}\right)^{-1}\\
0&\left(\begin{array}{cc}C&0\\ \gamma&C\end{array}\right)^{-1}
\end{array}\right)\left(\begin{array}{c}x\\ X\\ y\\ Y\end{array}\right).$$
Consequently, if we set $i:M\hookrightarrow M^c$ and $\pi:M^c\twoheadrightarrow M$, we immediately get  $\pi\circ\phi^c\circ i=\phi$.

\ \\




\end{document}